\input amstex
\magnification 1200
\TagsOnRight
\def\qed{\ifhmode\unskip\nobreak\fi\ifmmode\ifinner\else
 \hskip5pt\fi\fi\hbox{\hskip5pt\vrule width4pt
 height6pt depth1.5pt\hskip1pt}}
\NoBlackBoxes
\baselineskip 20 pt
\parskip 7 pt

\centerline {\bf RECONSTRUCTION OF THE WAVE SPEED}
\vskip -5 pt
\centerline {\bf FROM TRANSMISSION EIGENVALUES}
\vskip -5 pt
\centerline {\bf FOR THE SPHERICALLY-SYMMETRIC}
\vskip -5 pt
\centerline {\bf VARIABLE-SPEED WAVE EQUATION}

\vskip 7 pt
\centerline {Tuncay Aktosun}
\vskip -8 pt
\centerline {Department of Mathematics}
\vskip -8 pt
\centerline {University of Texas at Arlington}
\vskip -8 pt
\centerline {Arlington, TX 76019-0408, USA}
\vskip -8 pt
\centerline {aktosun\@uta.edu}

\centerline {Vassilis G. Papanicolaou}
\vskip -8 pt
\centerline {Department of Mathematics}
\vskip -8 pt
\centerline {National Technical University of Athens}
\vskip -8 pt
\centerline {Zografou Campus}
\vskip -8 pt
\centerline {157 80, Athens, Greece}
\vskip -8 pt
\centerline {papanico\@math.ntua.gr}

\noindent {\bf Abstract}: The unique reconstruction of a spherically-symmetric
wave speed $v$ is considered in a bounded spherical
region of radius $b$
 from the set of corresponding transmission
eigenvalues for which the corresponding eigenfunctions are
also spherically symmetric. If the integral of $1/v$ on the interval
$[0,b]$ is less than $b,$ assuming that there exists at least one
$v$ corresponding to the data, $v$ is uniquely
reconstructed from the data consisting of such transmission
eigenvalues and their ``multiplicities," where the
multiplicity is defined as the multiplicity of the transmission
eigenvalue as a zero of a key quantity. When that integral is
equal to $b,$ the unique reconstruction is presented when the data set
contains one additional piece of information. Some similar
results are presented for the unique reconstruction of the
potential from the transmission eigenvalues with
multiplicities for a related Schr\"odinger equation.

\vskip -8 pt
\par\noindent {\bf Short title:} Inverse problem for transmission eigenvalues
\vskip -8 pt
\par\noindent {\bf Keywords:}
transmission eigenvalues, inverse spectral problem,
variable-speed wave equation, Schr\"odinger equation, reconstruction, uniqueness

\newpage

\noindent {\bf 1. INTRODUCTION}
\vskip 3 pt

Let us consider the mathematical problem [8-10]
$$\cases
\Delta \Psi +\lambda\,\rho(\bold x)\,\Psi=0,\qquad \bold x\in\Omega,\\
\noalign{\medskip}
\Delta\Psi_0+\lambda\Psi_0=0,\qquad\bold x\in\Omega,\\
\noalign{\medskip} \Psi=\Psi_0,\quad \displaystyle\frac{\partial
\Psi}{\partial\bold n}=\frac{\partial\Psi_0}{\partial\bold
n}, \qquad \bold x\in \partial \Omega, \endcases \tag 1.1$$
where $\Delta $ denotes the Laplacian, $\lambda$ is the
spectral parameter, $\Omega$ is a bounded and simply connected domain
in ${\bold R}^n$ for any positive integer $n$
with the sufficiently smooth boundary $\partial \Omega,$
$\bold n$ is the outward unit
normal to the boundary $\partial \Omega,$ and the positive
quantity $\rho(\bold x)$ is assumed to be $1$ outside $\Omega.$
Those $\lambda$-values for which there are nontrivial solutions
$\Psi$ and $\Psi_0$ to (1.1) are known as the transmission
eigenvalues for (1.1). Obviously, $\lambda=0$ is always a transmission
eigenvalue, which we view as the trivial one.
The problem described in (1.1) arises in acoustic scattering
 from a bounded region with $1/\sqrt{\rho(\bold x)}$
denoting the wave speed $v$
and also
in electromagnetic
scattering from a bounded nonhomogeneity with refractive index
$\sqrt{\rho(\bold x)}$ as a function of location.

The relevant direct problem involves
the determination of the transmission eigenvalues
when the nonhomogeneity $\rho$ is known. The relevant inverse problem is to determine
the nonhomogeneity everywhere in the given domain $\Omega$ by using
an appropriate set of $\lambda$-values related to
the transmission eigenvalues
of (1.1).
It is already
known that the transmission eigenvalues for (1.1) can be
determined from some far-field measurements [4,5].
Letting $\lambda=k^2,$ we see that, corresponding to each transmission eigenvalue
$\lambda,$ we have two $k$-values, namely $k$ and $-k.$

The research field of direct and inverse problems involving
transmission
eigenvalues is now very active,
 and the related literature is growing rapidly, and hence
 it is impossible to provide a complete
bibliography on the general topic of transmission eigenvalues.
We refer the reader to [10] and the references therein to trace the important developments in the field and to [10,15-17]
and the references therein for the related
inverse problem of recovery of a bounded nonhomogeneity from an appropriate
set of transmission eigenvalues.

We are interested in (1.1) in the special case where $\Omega$ is a sphere of radius
$b$ centered at the origin in ${\bold R}^3$ and $\rho(\bold x)$
is spherically symmetric, which we write as $\rho(x)$ with $x:=|\bold x|.$
We further consider a subset of the transmission eigenvalues
for which the corresponding
eigenfunctions are also spherically symmetric. As in [2] we refer to
such $\lambda$-eigenvalues as {\it special transmission eigenvalues} of (1.1).

Let us assume that $\rho$ belongs to the {\it admissible
class} $\Cal A,$ by which we mean that
$\rho(x)$ for $x\in[0,+\infty)$ is positive, continuously differentiable,
and equal to $1$ for $x\ge b,$ and
that $\rho''(x)$ exists almost everywhere, where the prime denotes the $x$-derivative.

Let us introduce our key quantity
$$D(k):=\frac{\sin(kb)}{k}
\,\phi'(b;k)-\cos (kb)\,\phi(b;k),\tag 1.2$$
where $\phi(x;k)$ is the unique solution to
$$\cases\phi''+k^2\, \rho (x)\,\phi=0,\qquad 0<x<b,\\
\noalign{\medskip} \phi(0) = 0, \quad \phi'(0)=1.\endcases  \tag 1.3$$
Note that $D(-k)=D(k),$ and hence $D(k)$ is actually a function
of $k^2.$

The following result is known [2], and here we restate it in terms of
$k$ rather than $\lambda$ for later use.

\noindent {\bf Theorem 1.1} {\it Consider the
special case of (1.1) with
$\Omega$ being the three-dimensional ball of
radius $b$ centered at the origin,
where only
spherically-symmetric wavefunctions are allowed
and it is assumed that such wavefunctions
are continuous in the closure of $\Omega.$
Suppose that $\rho$ belongs to the admissible class $\Cal A.$ Then, the
corresponding special transmission eigenvalues of
(1.1) coincide with the $k^2$-values related to
the zeros of the quantity $D(k)$ defined in (1.2).}

As in [2], with each nonzero special transmission
eigenvalue $k_n^2$ of (1.1), we associate
a multiplicity, which is the same as the multiplicity of $k_n$ as a zero of $D(k).$
When $\rho$ belongs to the admissible class $\Cal A,$ it is known [2] that
$D(k)$ is entire in $k^2$ and has the representation
$$D(k)=\gamma\,E(k),\quad
E(k):= k^{2d}\prod_{n=1}^{\infty }\left( 1-\frac{
k^2}{k_n^2}\right),\tag 1.4$$
where $\gamma$ is a real constant, $k_n^2$ for $n\in\bold N$
correspond to the nonzero
transmission eigenvalues, and $d$ is the multiplicity
of the trivial zero transmission eigenvalue.
Here we use $\bold N:=\{1,2,3,\dots\}.$
It is known [2] that some $k_n^2$ may be repeated,
$d$ is at least $1,$ and the actual
value of $d$ is determined by $\rho.$ In the trivial case $\rho(x)\equiv 1,$ we
have $\gamma=0$ and hence $D(k)\equiv 0.$

In analyzing a typical transmission eigenvalue problem, one has to deal with a nonselfadjoint eigenvalue problem [2,8-10] and hence in general we cannot expect that all eigenvalues will be real. The nonselfadjointness in a natural way forces us to consider complex transmission eigenvalues and also transmission eigenvalues with multiplicities. The mathematical necessity of including complex
and nonsimple transmission eigenvalues certainly complicates the analysis. The experimental determination of such transmission eigenvalues, especially in the presence of nonsimple eigenvalues, presents a challenge and it is an interesting and important open question how to measure them and how to determine their multiplicities.

The constant $a$ defined as
$$a:=\int_0^b dx\,\sqrt{\rho (x)},\tag 1.5$$
has the physical interpretation as the travel
time for the wave to move from $x=0$ to $x=b.$
In [2] we have presented the following uniqueness results
regarding the determination of $\rho$ corresponding to
the special transmission eigenvalues of (1.1) in the spherically symmetric case.
Note that it is assumed that we know the value
of $b,$ which is a reasonable assumption as far as the applications are concerned.

\noindent {\bf Theorem 1.2} {\it
Consider the
special case of (1.1) with
$\Omega$ being the three-dimensional ball of
radius $b$ centered at the origin,
where only
spherically-symmetric wavefunctions are allowed
and it is assumed that such wavefunctions
are continuous in the closure of $\Omega.$
Suppose that our data set
consists of the corresponding special transmission
eigenvalues with their multiplicities, and assume that
there exists at least one corresponding
$\rho$ in the admissible class $\Cal A.$
Let $a$ be the constant defined in (1.5). We have the following:}

\item{(a)} {\it If $a<b,$ then our data set uniquely
determines $\rho;$ in other words, if both $\rho_1$ and $\rho_2$
correspond to our data, then we must have
$\rho_1\equiv \rho_2.$}

\item{(b)} {\it If $a=b,$ then our data set along with the value of
$\gamma$ appearing in (1.4) uniquely
determines $\rho;$ in other words, if both $\rho_1$ and $\rho_2$
correspond to our data and to the same
$\gamma,$ then we must have
$\rho_1\equiv \rho_2.$}

In this paper we give an alternate
proof of Theorem~1.2, a proof different from
that given in [2], by providing an algorithm
to reconstruct $\rho$ from the relevant data set.
Our paper is organized as follows. In Section~2 we present
some preliminary results that are needed later on;
the key result in Theorem~2.3 is crucial
for the unique reconstruction of $\rho$ given in Sections~3 and 4.
In Section~3 we present the alternate proof of Theorem~1.2(a) and
the reconstruction when $a<b.$ In Section~4 we present the alternate proof of Theorem~1.2(b) and
the reconstruction when $a=b.$
In Section~5 we consider the analogous problem for the Schr\"odinger
equation. We present Theorem~5.2, which is the analog
of Theorem~2.3 and which plays a key role
in the reconstruction of the potential in the Schr\"odinger equation. We then give an alternate proof of
the uniqueness result of Theorem~5.4,
a proof different from that given in [2], by providing
a reconstruction procedure for the potential
in terms of
the data set consisting of the corresponding transmission eigenvalues
with their multiplicities and the
parameter $\tilde\gamma$ appearing in (5.15).
In Section~5 we also provide an illustrative example showing that
we cannot have uniqueness if the data set does not include
$\tilde\gamma;$ let us mention, though, that the reconstructed potential
is outside the admissible class of potentials $\tilde A$
considered in Section~5 of our paper.
Finally, in Section~6 we present some explicit examples where
we display $\rho(x)$ and the corresponding quantities
$\gamma$ and $E(k)$ appearing in (1.4); however, in each of those examples,
$\rho'(x)$ has a jump discontinuity and hence is outside the admissible class
$\Cal A.$ In one of the examples presented in Section~6 it is shown that the same $E(k)$
yields two distinct $\rho(x)$ quantities, for one of which we have $a<b$ and for the other we have
$b>a.$

In the recovery algorithms given in Sections~3-5, we solve some basic Riemann-Hilbert problems
and use basic facts related to their unique solutions. For the benefit
of the readers who are unfamiliar with the theory of
Riemann-Hilbert problems, we summarize below their formulation and their unique
solutions relevant to our paper.
Let us use ${\bold C}$ for the complex plane, ${\bold C^+}$ for the open upper-half
complex plane, $\overline{\bold C^+}$ for ${\bold C^+}\cup{\bold R},$ ${\bold C^-}$ for the open lower-half
complex plane, and $\overline{\bold C^-}$ for ${\bold C^-}\cup{\bold R}.$

The idea behind solving a basic Riemann-Hilbert problem is to determine
a sectionally analytic function on ${\bold C}$ by determining
its sections on ${\bold C^+}$ and on ${\bold C^-},$ respectively, from its jump value
on the real axis ${\bold R}.$ Mathematically, we need to solve the
functional equation
$$F(k)-F(-k)=G(k),\qquad k\in{\bold R},\tag 1.6$$
where $G(k)$ is relevant only for real values of $k$
and it indicates the jump. In other words, given $G(k)$ for
$k\in{\bold R},$ we need to determine $F(k)$ for $k\in{\bold C^+}$ in such a way that
$F(k)$ is analytic in $k\in{\bold C^+},$ continuous in $k\in\overline{\bold C^+},$ and
$O(1/k)$ as $k\to\infty$ in $\overline{\bold C^+}.$ Certainly, we then have
$F(-k)$ analytic in $k\in{\bold C^-},$ continuous in $k\in\overline{\bold C^-},$ and
$O(1/k)$ as $k\to\infty$ in $\overline{\bold C^-}.$
In general $G(k)$ may not
have any extension off the real axis, but even if it does
only the values of $G(k)$ for $k\in{\bold R}$ are relevant and needed
in solving (1.6). For the unique solvability of (1.6) it is sufficient
to assume that $G(k)$ behaves like $O(1/k)$ as $k\to\pm\infty$
on the real axis ${\bold R}$ and that $G(k)$ is H\"older continuous on ${\bold R}$
with a positive index $\alpha.$ The latter condition is
expressed as
$$|G(k_1)-G(k_2)|\le C|k_1-k_2|^\alpha,\qquad k_1,k_2\in{\bold R},\tag 1.7$$
for some positive constant $C$ independent of $k.$ In the special case $\alpha=1$ the
condition given in (1.7) is known as the Lipschitz continuity
of $G(k)$ on ${\bold R}.$ In fact, in our paper the relevant
$G(k)$ is Lipschitz continuous on ${\bold R}$ due to the fact that
$G(k)$ has an analytic continuation from $k\in{\bold R}$ to
$k\in{\bold C};$ even though the relevant $G(k)$ is bounded on
${\bold R}$ and decays as $O(1/k)$ as $k\to\pm\infty$ on ${\bold R},$
in general $G(k)$ grows exponentially as $k\to\infty$ in ${\bold C^+}$
and in ${\bold C^-},$ but as already stated that unboundedness
off the real axis is irrelevant in the analysis
of the Riemann-Hilbert problem given in (1.6).
Under the two aforementioned sufficiency conditions on
$G(k)$ for $k\in{\bold R},$ the Riemann-Hilbert problem
given in (1.6) is uniquely solvable and the analytic section
$F(k)$ defined on $k\in\overline{\bold C^+}$ is explicitly expressed in terms of
the values of $G(k)$ known for $k\in{\bold R}$ as
$$F(k)=\displaystyle\frac{1}{2\pi i}\displaystyle\int_{-\infty}^\infty dt\,\displaystyle\frac{G(t)}{t-k-i0^+},
\qquad k\in\overline{\bold C^+},\tag 1.8$$
where one can ignore $i0^+$ in (1.7) if $k\in{\bold C^+}$ and one can
interpret the presence of $i0^+$ by stating that
$F(k)$ for $k\in{\bold R}$ must be obtained by first evaluating
the integral in (1.8) when $k\in{\bold C^+}$ and then letting $k$ approach
its real value from ${\bold C^+}.$ When
$G(k)$ satisfies the two relevant conditions on ${\bold R},$
it may sometimes be possible to solve
(1.6) readily by finding a familiar function $F(k)$ satisfying
(1.6) with the appropriate properties on $\overline{\bold C^+}.$ Since
the existence and uniqueness are ensured, we can then conclude that
that function $F(k)$ must satisfy (1.8).
We refer the reader to [11,13] for
further information on Riemann-Hilbert problems and their solutions.

In our paper the Riemann-Hilbert problem stated in (1.6) arises in
(3.4), (4.5), and (5.14). For example, the basic idea behind solving
the Riemann-Hilbert problem in (4.5) is to split $2ikD(k)$ into the two pieces
$Q(k)$ and $Q(-k),$ where we know that $2ikD(k)$ has the behavior $O(1/k)$ as $k\to\pm\infty$ on ${\bold R}$
and is Lipschitz continuous
for $k\in{\bold R}.$
In our specific case it turns out that $2ikD(k)$ has an analytic extension in $k$ to the entire
complex plane ${\bold C}$ with an exponential growth as
$k\to\infty$ in ${\bold C}$ except when $k\to\pm\infty$
on the real axis. The splitting is such that $Q(k)$ is analytic in $k\in{\bold C^+}$ and
bounded in $\overline{\bold C^+}$ and in fact $Q(k)=O(1/k)$ as $k\to\infty$ in $\overline{\bold C^+}.$
Similarly, the other piece $Q(-k)$ is analytic in $k\in{\bold C^-}$ and
bounded in $\overline{\bold C^-}$ and in fact $Q(-k)=O(1/k)$ as $k\to\infty$ in $\overline{\bold C^-}.$

In Sections~3 and 4, the cases $a<b$ and $a=b,$ respectively,
are analyzed, where $a$ and $b$ are the constants appearing in (1.5).
In both cases we show that the relevant $G(k)$ in (1.6) can be split
into the appropriate functions
$F(k)$ and $F(-k)$ in such a way that
we are able to recognize $F(k)$ and express it explicitly
in terms of a familiar spectral function, and hence
we are able to
reconstruct the bounded nonhomogeneity from
the explicitly constructed $F(k).$
On the other hand, when $a>b,$ even though
we know that the relevant Riemann-Hilbert given in (1.6)
is uniquely solvable, we are unable to express
the corresponding $F(k)$ explicitly in terms of a familiar
spectral function yielding the nonhomogeneity. Thus, our reconstruction method
exploiting the relevant Riemann-Hilbert given in (1.6) does not seem to
yield the nonhomogeneity when $a>b.$ In other words, the case $a>b$ is an open problem,
and it is not known in that case whether the nonhomogeneity can be recovered
by a method similar to that used in Sections~3 and 4.

\vskip 8 pt
\noindent {\bf 2. PRELIMINARIES}
\vskip 3 pt

We consider the extension of the differential equation in
(1.3) to the half line
${\bold R}^+:=(0,+\infty),$ namely
$$\psi''+k^2\rho(x)\,\psi=0,\qquad x\in{\bold R}^+,\tag 2.1$$
where
$\rho(x)$ belongs to the admissible class $\Cal A$ and hence
$\rho(x)\equiv 1$ for $x\ge b.$ Let us define
the travel-time coordinate $y$ as
$$y(x):=\int_0^x ds\,\sqrt{\rho(s)},\qquad
x\in[0,+\infty),\tag 2.2$$
and note that (1.5) and (2.2) imply that
$a=y(b).$ We remark that (2.2) yields
$$y(x)=x+a-b,\qquad x\ge b.$$

Let $f(x;k)$ denote the Jost solution to (2.1), i.e. $f(x;k)$ satisfies (2.1) and
$$f(x;k)=e^{ikx},\quad
f'(x;k)=ik\,e^{ikx},\qquad x\ge b.\tag 2.3$$
We also note that, when $k=0,$
(2.1) reduces to $\psi''(x;0)=0,$ and hence with the help of (2.3) we get
$$f(x;0)=1,\qquad x\in[0,+\infty).\tag 2.4$$

Via a Liouville transformation, (2.1) can be transformed into
a Schr\"odinger equation. In other words, if we let
$$\tilde f(y;k):=[\rho(x)]^{1/4}\,e^{-ik(b-a)}\,f(x;k),\tag 2.5$$
where $y$ is related to $x$ as in (2.2),
then $\tilde f(y;k)$ becomes the Jost solution to the
Schr\"odinger equation
$$\tilde f''(y;k)+k^2 \tilde f(y;k)=V(y)\,\tilde f(y;k),\qquad y\in{\bold R}^+,\tag 2.6$$
where the prime now denotes the $y$-derivative and we have
$$\tilde f(y;k)=e^{iky},\quad
\tilde f'(y;k)=ik\,e^{iky},\qquad y\ge a,\tag 2.7$$
$$V(y):=V(y(x))=\displaystyle\frac{\rho''(x)}{4\,[\rho(x)]^2}-\displaystyle\frac{5 [\rho'(x)]^2}{16\,[\rho(x)]^3} .\tag 2.8$$

With the help of (2.2) and a multiple use of the chain rule
in taking the derivatives
on the right-hand side of (2.8), we write (2.8) as the second-order
linear differential equation
$$\displaystyle\frac{d^2 \left(\left[\rho(x(y))\right]^{1/4}\right)}{dy^2}=V(y)\,[\rho(x(y))]^{1/4},
\qquad y\in{\bold R}^+.
\tag 2.9$$
Note that $\rho(x(y))$ satisfies
$$\rho(x(y))\big|_{y=a}=1,\quad
\displaystyle\frac{d[\rho(x(y))]}{dy}\bigg|_{y=a}=0,$$
and hence
$$[\rho(x(y))]^{1/4}\big|_{y=a}=1,\quad
\displaystyle\frac{d[\rho(x(y))]^{1/4}}{dy}\bigg|_{y=a}=0.$$
As we shall see,
(2.9) will be useful
in constructing explicit illustrative examples of $\rho(x)$ and $V(y)$
 from some appropriate sets of data.

 From (2.4) and (2.5) we conclude the following result.

\noindent {\bf Corollary 2.1} {\it Assume that
$\rho(x)$ belongs to the admissible class $\Cal A.$ Let
$\tilde f(y;k)$ be the corresponding Jost solution to (2.6)
with $V(y)$ as in (2.8).
Then, we have}
$$\tilde f(y;0)=[\rho(x(y))]^{1/4},\qquad y\in[0,+\infty).\tag 2.10$$

Let us remark that, when $\rho(x)$ is in the admissible class $\Cal A,$ we have [6,12,14,18]
$$f(x;-k)=f(x;k)^*,\qquad k\in{\bold R},\tag 2.11$$
where the asterisk denotes complex conjugation.
The result in (2.11) will be useful in establishing (2.18) and
in the reconstruction of $\rho(x)$ from a data set
containing the associated transmission eigenvalues.

The results in the following proposition are already known but we
state them with a brief proof for the convenience of the reader.

\noindent {\bf Proposition 2.2} {\it
Assume that
$\rho(x)$ belongs to the admissible class $\Cal A.$ Let
$V(y)$ be the potential obtained from $\rho(x)$ as in (2.8)
and let
$\tilde f(y;k)$ be the corresponding Jost solution to (2.6) satisfying (2.7).
Then, we have the following:}

\item{(a)} {\it The potential $V(y)$ belongs to the admissible class
$\tilde A$ described in Section~5. Consequently, the Jost solution
$\tilde f(y;k)$ has the properties outlined in Proposition~5.1.}

\item{(b)} {\it The differential equation (2.1) with the
Dirichlet boundary condition $\psi(0)=0,$ cannot have,
for any negative value of $k^2,$ any
solutions that are square integrable in $x\in{\bold R}^+.$}

\item{(c)} {\it The corresponding half-line Schr\"odinger equation (2.6)
with the Dirichlet boundary condition at $y=0$ cannot have any nontrivial
solutions that are square integrable in $y\in{\bold R}^+.$
Hence, the corresponding Schr\"odinger operator has no bound states,
and therefore
 $\tilde f(0;k)$ cannot vanish for $k\in\overline{\bold C^+}\setminus\{0\}.$}

\noindent PROOF: When $\rho$ is in the admissible class $\Cal A,$
the corresponding $V(y),$ as seen from (2.2) and (2.8),
is real valued, compactly supported, and
integrable. Thus, $V(y)$ belongs to the admissible class
$\tilde A$ described in Section~5, and (a) is proved. Note that (c)
directly follows from (b) because of the Liouville transformation given
in (2.5), Proposition~5.1, and the fact that the bound states correspond to
square-integrable solutions of the relevant
differential equations. Hence, we only need to prove (b).
Because of the selfadjointness of the corresponding Schr\"odinger
operator, any existing bound states may occur only at negative
values of $k^2.$
If $\psi(x)$ were a nontrivial square-integrable solution to
(2.1) at some negative value of
$k^2,$ i.e. if $\psi\in L^2({\bold R}^+),$ then from (2.1) we would get
$\psi''\in L^2({\bold R}^+),$ and hence,
e.g. via Fourier transforms, we would have
$\psi'\in L^2({\bold R}^+).$ By the Cauchy-Schwarz inequality, we would then have
$\psi \psi'\in L^1({\bold R}^+).$ However, that would imply the existence of a
sequence $x_n$ converging to $+\infty$ such that
$$|\psi'(x_n)\,\psi(x_n)|\to 0,\qquad x_n\to+\infty.$$
 From (2.1) through integration we would then get
$$\int_0^{x_n} dx\, [\psi(x)]^*\psi''(x)+k^2
\int_0^  {x_n} dx\,\rho(x)\,|\psi(x)|^2=0.\tag 2.12$$
Using the Dirichlet condition $\psi(0)=0$ and an
integration by parts on the first integral in
(2.12),  we would obtain
$$[\psi(x_n)]^*\psi'(x_n)-[\psi(0)]^*\psi'(0^+)
-\displaystyle\int_0^{x_n} dx\,|\psi'(x)|^2
+k^2
\displaystyle\int_0^{x_n} dx\,\rho(x)\,|\psi(x)|^2=0.\tag 2.13$$
Since $\psi(0)=0,$ from (2.1) it follows that $\psi''(0^+)=0$ and hence
 $\psi'(0^+)$ is finite. Thus, letting
 $x_n\to+\infty,$ from (2.13) we would get
$$-\displaystyle\int_0^{\infty} dx\,|\psi'(x)|^2
+k^2
\displaystyle\int_0^{\infty} dx\,\rho(x)\,|\psi(x)|^2=0,\tag 2.14$$
which is a contradiction
because the left-hand side of (2.14) would be strictly negative
due to the fact that $k^2<0,$
$\psi(x)$ is assumed to be
a nontrivial solution, and $\rho(x)>0$ for $x\in{\bold R}^+.$ \qed

Let $\phi(x;k)$ and $\tilde \phi(y;k)$ denote the solutions to the
initial-value problems on the half line that are respectively given by
$$\cases \phi''(x;k)+k^2\rho(x)\,\phi(x;k)=0,\qquad x\in{\bold R}^+,\\
\noalign{\medskip}
\phi(0;k)=0,\quad \phi'(0;k)=1,\endcases\tag 2.15$$
$$\cases \tilde\phi''(y;k)+k^2\tilde\phi(y;k)=V(y)\,\tilde \phi(y;k),\qquad y\in{\bold R}^+,\\
\noalign{\medskip}
\tilde\phi(0;k)=0,\quad \tilde\phi'(0;k)=1,\endcases\tag 2.16$$
where $V(y)$ is related to $\rho(x)$ as in (2.8) and $x$ and $y$ are
related to each other as in (2.2).
Note that (2.15) and (2.16) are uniquely solvable [7] and that the corresponding solutions
are entire in $k^2.$
We remark that (2.15) is actually the extension of (1.3) from the interval $(0,b)$
to ${\bold R}^+$ and hence we use $\phi(x;k)$ to denote the unique solution
to both (1.3) and (2.15).

The result in the following theorem is crucial for the reconstruction of $\rho(x)$
 from the data containing transmission eigenvalues, and it will be used in Sections~3 and 4.

\noindent {\bf Theorem 2.3} {\it Let $\rho(x)$ belongs to the
admissible class $\Cal A.$ Then,
the quantity $D(k)$ defined in (1.2) is related to the Jost solution $f(x;k)$
to (2.1) as}
$$D(k)=\displaystyle\frac{f(0;k)-f(0;-k)}{2ik},\qquad k\in{\bold C}.\tag 2.17$$
{\it For real $k$-values, we
have}
$$\text{\rm{Im}}[f(0;k)]=k\,D(k),\qquad k\in{\bold R},\tag 2.18$$
{\it where $\text{\rm{Im}}$ denotes the imaginary part.}

\noindent PROOF:
Let us express the solution $\phi(x;k)$ to (2.15) as a linear
combination of the
linearly independent solutions $f(x;k)$ and $f(x;-k)$ to (2.1),
where $f(x;k)$ is the Jost solution satisfying (2.3).
We have
$$\bmatrix \phi(x;k)\\
\noalign{\medskip}
\phi'(x;k)\endbmatrix=\bmatrix f(x;k)&f(x;-k)\\
\noalign{\medskip}
f'(x;k)&f'(x;-k)\endbmatrix \bmatrix c_1(k)\\
\noalign{\medskip}
c_2(k)\endbmatrix,\tag 2.19$$
where the coefficients $c_1(k)$ and $c_2(k)$ are independent of $x$
and are yet to be determined.
With the help of (2.3) and the second line of (2.15), we evaluate (2.19)
at $x=b$ and $x=0,$ respectively, and we obtain
$$\bmatrix \phi(b;k)\\
\noalign{\medskip}
\phi'(b;k)\endbmatrix=\bmatrix e^{ikb}&e^{-ikb}\\
\noalign{\medskip}
ik\,e^{ikb}&-ik\,e^{-ikb}\endbmatrix \bmatrix c_1(k)\\
\noalign{\medskip}
c_2(k)\endbmatrix,\tag 2.20$$
$$\bmatrix 0\\
1\endbmatrix=\bmatrix f(0;k)&f(0;-k)\\
\noalign{\medskip}
f'(0;k)&f'(0;-k)\endbmatrix \bmatrix c_1(k)\\
\noalign{\medskip}
c_2(k)\endbmatrix.\tag 2.21$$
 From (2.20) and (2.21), by eliminating $c_1(k)$ and $c_2(k)$ we get
$$\bmatrix \phi(b;k)\\
\noalign{\medskip}
\phi'(b;k)\endbmatrix=\bmatrix e^{ikb}&e^{-ikb}\\
\noalign{\medskip}
ik\,e^{ikb}&-ik\,e^{-ikb}\endbmatrix
\bmatrix f(0;k)&f(0;-k)\\
\noalign{\medskip}
f'(0;k)&f'(0;-k)\endbmatrix^{-1}
\bmatrix 0\\
\noalign{\medskip}
1\endbmatrix.\tag 2.22$$
Let $[g;h]:=gh'-g'h$ denote the Wronskian of any two functions
$g$ and $h.$ It
is known [7] and can also directly be verified that the Wronskian of any
two solutions to (2.1) is independent of $x.$ With the help of
(2.1) and (2.3) we get
$$[f(x;k);f(x;-k)]=-2ik,$$
and hence
$$\bmatrix f(0;k)&f(0;-k)\\
\noalign{\medskip}
f'(0;k)&f'(0;-k)\endbmatrix^{-1}=\displaystyle\frac{-1}{2ik}
\bmatrix f'(0;-k)&-f(0;-k)\\
\noalign{\medskip}
-f'(0;k)&f(0;k)\endbmatrix.\tag 2.23$$
Using (2.23) in (2.22) we obtain
$$\bmatrix \phi(b;k)\\
\noalign{\medskip}
\phi'(b;k)\endbmatrix=\displaystyle\frac{-1}{2ik}\,
\bmatrix e^{ikb}&e^{-ikb}\\
\noalign{\medskip}
ik\,e^{ikb}&-ik\,e^{-ikb}\endbmatrix
\bmatrix -f(0;-k)\\
\noalign{\medskip}
f(0;k)\endbmatrix.\tag 2.24$$
Writing (1.2) as the matrix product
$$D(k)=\bmatrix -\cos(kb)&\displaystyle\frac{\sin(kb)}{k}\endbmatrix
\bmatrix \phi(b;k)\\
\noalign{\medskip}
\phi'(b;k)\endbmatrix,\tag 2.25$$
and using (2.24) in (2.25), after some simplification we obtain
(2.17). Finally, using (2.11) in (2.17), we get (2.18). \qed

\vskip 8 pt
\noindent {\bf 3. RECONSTRUCTION OF $\rho(x)$ WHEN $a<b$}
\vskip 3 pt

In this section we give a
proof of Theorem~1.2(a) by providing a reconstruction algorithm for
the unique recovery of $\rho(x)$ in terms of the data consisting of the
corresponding special
transmission eigenvalues with their multiplicities. Thus, our
data set is equivalent to the set of zeros (including
the multiplicities of those zeros) of
the quantity $E(k)$ defined in (1.4).
Equivalently, the knowledge of our data is equivalent to
the knowledge of $E(k).$
Using (2.5) in (2.17), with the help of (1.4), we get
$$E(k)=
\displaystyle\frac{e^{ik(b-a)}\tilde f(0;k)-e^{-ik(b-a)}\tilde f(0;-k)}
{2ik\,\gamma\,[\rho(0)]^{1/4}}.\tag 3.1$$

Note that we assume that
$a<b,$ where $a$ is the constant defined in (1.5).
For the reconstruction, we assume that the existence problem
is solved, i.e. we assume the existence of at least
one $\rho$ in the admissible class $\Cal A$ corresponding to our data.
The uniqueness aspect in the recovery of
$\rho(x)$ follows from the uniqueness in each of the
reconstruction steps outlined below:

\item {(a)} When $\rho(x)$ is in the admissible class $\Cal A,$
as stated in Theorem~2.2(a), the
corresponding Jost solution $\tilde f(y;k)$ given in (2.5) satisfies
the properties listed in
Proposition~5.1, and in particular (5.3) holds.
Thus, using (5.3) in (3.1)
we get the large-$k$ asymptotics of
$E(k)$ for $k\in{\bold R}$ as
$$E(k)=
\displaystyle\frac{\sin k(b-a)}{k\,\gamma\,[\rho(0)]^{1/4}}+O\left(\displaystyle\frac{1}{k^2}\right),
\qquad k\to\pm\infty.\tag 3.2$$
Since $(b-a)$ is assumed to be positive, from (3.2)
we can determine the value of $(b-a)$
and hence $a$ and also
the value of $\gamma\,[\rho(0)]^{1/4}.$  Note that we do not
have the values of $\gamma$ and $[\rho(0)]^{1/4}$ separately, but as we will see
this
does not create an obstacle for the reconstruction of $\rho(x).$

\item {(b)} Letting
$$P(k):=\tilde f(0;k)-1,\tag 3.3$$
we write (3.1) as
$$e^{ik(b-a)}P(k)-e^{-ik(b-a)}P(-k)=\varphi(k),\qquad k\in{\bold R},\tag 3.4$$
where we have defined
$$\varphi(k):=2i\left[
k\,E(k)\,\gamma\,[\rho(0)]^{1/4}-\sin k(b-a)\right].\tag 3.5$$
By the previous step given in (a) above, we know that our data set
uniquely determines the value
of $\varphi(k).$ Furthermore, using (3.2) in (3.5) we get
$\varphi(k)=O(1/k)$ as $k\to\pm\infty.$
When $\rho$ is in the admissible class $\Cal A,$
it is known that $E(k)$ is entire in $k^2$ and hence also in $k.$
The Lipschitz continuity of $\varphi(k)$ for $k\in{\bold R}$ follows from the
fact that the right hand side in (3.5) has an analytic extension
to the entire complex plane and that
$\varphi(k)=O(1/k)$ as $k\to\pm\infty.$

\item {(c)} Note that (3.4) constitutes a Riemann-Hilbert
problem on the complex plane where the function $\varphi(k)$
is specified for $k\in{\bold R}$ and it satisfies the Lipschitz
continuity on ${\bold R}$ and behaves as $O(1/k)$ as $k\to\pm\infty$
on ${\bold R}.$
The goal is to obtain
$e^{ik(b-a)}P(k)$ and $e^{-ik(b-a)}P(-k)$
in such a way that $e^{ik(b-a)}P(k)$ is analytic in ${\bold C^+},$
continuous in $\overline{\bold C^+},$ and $O(1/k)$ as $k\to\infty$
in $\overline{\bold C^+}.$ Since $b>a,$
those properties of $e^{ik(b-a)}P(k)$ follow
if and only if $P(k)$ satisfies those properties.
The unique solvability of (3.4) follows from the
Lipschitz continuity of $\varphi(k)$ on ${\bold R}$ and
the fact that $\varphi(k)=O(1/k)$ as $k\to\pm\infty.$
When $\rho$ in the admissible class, we already know
 from Proposition~5.1 that $\tilde f(0;k)$ is analytic in $k\in{\bold C^+},$ is
continuous in $k\in\overline{\bold C^+},$ and satisfies (5.3). Thus, the function
$P(k)$ given in (3.3) helps us to obtain the unique solution to
(3.4).
As we have indicated in Section~1, the unique solution to the
Riemann-Hilbert problem given in (3.4) is then expressed
with the help of (1.8) as
$$e^{ik(b-a)}P(k)=\displaystyle\frac{1}{2\pi i}\int_{-\infty}^\infty dt\,\displaystyle\frac{\varphi(t)}{t-k-i0^+},\qquad
k\in \overline{\bold C^+}.\tag 3.6$$
We see from (3.3) and (3.6) that
$$\tilde f(0;k)=1+\displaystyle\frac{e^{-ik(b-a)}}{2\pi i}\int_{-\infty}^\infty dt\,\displaystyle\frac{\varphi(t)}{t-k-i0^+},\qquad
k\in \overline{\bold C^+}.$$

\item {(d)} Having constructed $\tilde f(0;k)$ from $E(k),$ we
remark that, as stated in Theorem~2.2(c) $\tilde f(0;k)$
does not have any zeros on the positive imaginary axis
in the complex $k$-plane and hence the
corresponding half-line Schr\"odinger equation with the
Dirichlet boundary condition has no bound states.
Thus, we can use the Marchenko
procedure described in Proposition~5.1 to
uniquely reconstruct $V(y)$ of (2.8) and the Jost
solution $\tilde f(y;k)$ from the already
reconstructed $\tilde f(0;k).$
This is done by first
constructing the scattering matrix $\tilde S(k)$ from
$\tilde f(0;k)$ as in (5.5). Next, the Marchenko kernel
$M(\xi)$ is constructed as in (5.9) but without the summation
term there due to the fact that there are no bound states.
Then, the Marchenko integral equation (5.8) is uniquely
solved for $K(y,\xi),$ and the quantities
$V(y)$ and $\tilde f(y;k)$ are uniquely recovered as
in (5.10) and (5.11), respectively.

\item {(e)} Having recovered $\tilde f(y;k),$ we obtain
$\rho(x(y))$ from $\tilde f(y;0)$ with the help of (2.10),
namely
$$\rho(x(y))=[\tilde f(y;0)]^4,\qquad y\in[0,+\infty).\tag 3.7$$
Our next task is to obtain $\rho$ in terms of $x$ by establishing
the relationship between $x$ and $y.$ From (2.2) we have
$$\displaystyle\frac{dy}{dx}=\sqrt{\rho(x(y))}.\tag 3.8$$
Hence, from (2.2) and (3.8)
we obtain the first-order, separable ordinary
differential equation
$$\displaystyle\frac{dy}{[\tilde f(y;0)]^2}=dx,\qquad x\in{\bold R}^+,\tag 3.9$$
with the initial condition $y(0)=0.$ By integrating (3.9), the relationship
between $x$ and $y$ is obtained as
$$x=\int_0^y \displaystyle\frac{ds}{[\tilde f(s;0)]^2},\qquad y\in[0,+\infty)
.\tag 3.10$$
Since $\rho(x)$ is assumed to be positive,
 from (3.10) it follows that the mapping
$x\mapsto y$ is one-to-one and onto on ${\bold R}^+.$ Having $x$ as a
function of $y$ in (3.10), we can invert it to get $y$ as a
function of $x.$ Thus, by using (3.10) in (3.7) we recover
$\rho(x)$ in terms of $x$ as
$$\rho(x)=[\tilde f(y(x),0)]^4,\qquad x\in[0,+\infty).$$
Thus, the reconstruction of
$\rho(x)$ for $x\in[0,+\infty)$ from $E(k)$ for $k\in{\bold R}$
is accomplished.

Finally, let us note that our procedure yields
the value of the constant $\gamma$ from the knowledge of $E(k).$ This is because
we already have the value of $\gamma\,[\rho(0)]^{1/4}$ from (3.2)
and we have the value of $\rho(0)$ from (3.7) evaluated at $y=0.$

\vskip 8 pt
\noindent {\bf 4. RECONSTRUCTION OF $\rho(x)$ WHEN $a=b$}
\vskip 3 pt

In this section we consider the case $a=b,$
where $a$ and $b$ are the quantities appearing in (1.5).
We give an independent
proof of Theorem~1.2(b) by providing a reconstruction algorithm for
$\rho(x)$ from the data consisting of the
corresponding special
transmission eigenvalues with their multiplicities and the
constant $\gamma$ appearing in (1.4). By Theorem~1.1 the knowledge of our
data is equivalent to knowing the zeros (with multiplicities) of the
quantity $D(k)$ given in (1.4) as well as the value of $\gamma$ there.
Hence, our data set is equivalent to the knowledge of $D(k).$

As seen from (1.4) and (3.2), if $D(k)=O(1/k^2)$ as $k\to\pm\infty,$ we can deduce that $a=b.$

When $a=b,$ let us outline the unique recovery of $\rho(x)$ from $D(k)$
given in (1.4).

\item {(a)} Since $a=b,$ from (1.4) and (3.1) we see that
$$2ik D(k)=
\displaystyle\frac{\tilde f(0;k)-\tilde f(0;-k)}
{[\rho(0)]^{1/4}},\qquad k\in{\bold C}.\tag 4.1$$
On the other hand, from (2.10) we know that
$$\tilde f(0;0)=[\rho(0)]^{1/4},\tag 4.2$$
and hence we rewrite (4.1) as
$$2ik D(k)=
\displaystyle\frac{\tilde f(0;k)-\tilde f(0;-k)}
{\tilde f(0;0)},\qquad k\in{\bold C}.\tag 4.3$$
We know from (4.2) that $\tilde f(0;0)$ is real and in fact
positive because $\rho(0)>0.$
When $\rho(x)$ is in the admissible class $\Cal A,$
by using (5.3) and (5.4) in (4.1) we conclude that
$2ik D(k)=O(1/k)$ as $k\to\pm\infty$ on ${\bold R}.$
Furthermore, we know that $D(k)$ has an analytic
extension to the entire complex plane ${\bold C}$ and hence
$2ik D(k)$ is Lipschitz continuous on ${\bold R}.$
As seen from (5.3) and (5.4),
$2ik D(k)$ is unbounded as $k\to\infty$ in
${\bold C^+}$ and in ${\bold C^-},$ but for the analysis of the
Riemann-Hilbert problem to be studied we need the
large-$k$ asymptotics of $2ik D(k)$ only on the real axis.

\item {(b)} Letting
$$Q(k):=\displaystyle\frac{\tilde f(0;k)-1}
{\tilde f(0;0)},\tag 4.4$$
we write (4.3) for real $k$-values as
$$Q(k)-Q(-k)=2ik D(k),\qquad k\in{\bold R}.\tag 4.5$$
Note that (4.5) constitutes a Riemann-Hilbert
problem on the complex plane where the function $2ik D(k)$
is specified for $k\in{\bold R}$ and it is
Lipschitz continuous on ${\bold R}$ and
behaves as $O(1/k)$ as $k\to\pm\infty$ on ${\bold R}.$
The goal is to obtain
$Q(k)$ and $Q(-k)$
in such a way that $Q(k)$ is analytic in ${\bold C^+},$ continuous in $\overline{\bold C^+},$ and
$O(1/k)$ as $k\to\infty$ in $\overline{\bold C^+}.$ Since the unique
solvability
of (4.5) is assured by the two relevant properties
of $2ik D(k)$ stated in (a), we know that
the function given in (4.4) must be that unique solution.
The corresponding properties of
$Q(k)$ are deduced from (4.4) by using the relevant properties
of $\tilde f(0;k)$ given in Proposition~5.1.
In particular, using Proposition~5.1(a) we establish the analyticity of
$Q(k)$ in ${\bold C^+};$ using in (4.4) Proposition~5.1(a) and the fact that
$\tilde f(0;0)>0$ we conclude the continuity of
$Q(k)$ in $\overline{\bold C^+};$ using (5.3) in (4.4) we
obtain $Q(k)=O(1/k)$ as $k\to\infty$ in $\overline{\bold C^+}.$
Therefore, as indicated in (1.8), the unique solution to the
Riemann-Hilbert problem given in (4.5) satisfies
$$Q(k)=\displaystyle\frac{1}{2\pi i}\int_{-\infty}^\infty dt\,\displaystyle\frac{2it\, D(t)}{t-k-i0^+},\qquad
k\in \overline{\bold C^+}.\tag 4.6$$
We see from (4.4) and (4.6) that
$$\tilde f(0;k)=1+\displaystyle\frac{\tilde f(0;0)}{\pi}\int_{-\infty}^\infty dt\,\displaystyle\frac{t\, D(t)}{t-k-i0^+},\qquad
k\in \overline{\bold C^+}.\tag 4.7$$

\item {(c)} Let us remark that (4.7) also follows from the
Schwarz integral formula for the half plane [1] by using the following argument. As stated earlier,
$Q(k)$ is analytic in ${\bold C^+},$ continuous in $\overline{\bold C^+},$ and
$O(1/k)$ as $k\to\infty$ in $\overline{\bold C^+}.$
For real $k$-values, from (4.5) and
(5.2), we obtain
$$k\,D(k)=\text{Im}[Q(k)],\qquad k\in{\bold R}.\tag 4.8$$
Thus, we can construct
$Q(k)$ for $k\in\overline{\bold C^+}$ from its imaginary part known for $k\in{\bold R}$
by using
the Schwarz integral formula [1]
$$Q(k)=\displaystyle\frac{1}{\pi}\int_{-\infty}^\infty dt\,\displaystyle\frac{\text{Im}[Q(t)]}{t-k-i0^+},\qquad
k\in \overline{\bold C^+}.\tag 4.9$$
Hence, (4.8) and (4.9) yield (4.6) and in turn (4.7).

\item {(d)} Having constructed $Q(k)$ from $D(k),$ we can evaluate
the value of $\tilde f(0;0)$ by setting $k=0$ in
(4.4), which yields
$$\tilde f(0;0)=
\displaystyle\frac{1}{1-Q(0)}.
\tag 4.10$$
Thus, (4.6), (4.7), and (4.10) imply that $D(k)$ uniquely determines
$\tilde f(0;k).$

\item {(e)} Having reconstructed $\tilde f(0;k)$ from $D(k),$
we can then repeat the reconstruction
steps (d) and (e) of Section~3 in order to uniquely
reconstruct $\rho(x)$. Thus, the reconstruction of $\rho(x)$ for
$x\in[0,b]$ with $b=a$ from the data set $D(k)$
is accomplished.

As stated in [2], when $b=a$ it is an open problem whether
the value of the constant $\gamma$ appearing in (1.4) is needed or whether $\gamma$ can be determined from the data set consisting of $E(k)$ given in (1.4) alone.

\vskip 8 pt
\noindent {\bf 5. RECONSTRUCTION OF THE POTENTIAL IN THE
SCHR\"ODINGER EQUATION}
\vskip 3 pt

The transmission eigenvalue problem for (1.1) has an analog
for the Schr\"odinger equation. The transmission eigenvalues in that case
correspond to those values of $k^2$
for which there exist a nontrivial solution pair $\Psi$ and
$\Psi_0$ to the system
$$
\cases
-\Delta \tilde\Psi + V(\bold y) \,\tilde\Psi=k^2\tilde\Psi,\qquad \bold y\in \Omega,\\
\noalign{\medskip}
-\Delta \tilde\Psi_0 = k^2 \tilde\Psi_0,\qquad \bold y\in \Omega, \\
\noalign{\medskip}
\tilde\Psi=\tilde\Psi_0,\quad
\displaystyle\frac{\partial \tilde\Psi}{\partial \bold n}
=\frac{\partial \tilde\Psi_0}{\partial \bold n},
\qquad \bold y\in \partial \Omega, \endcases \tag 5.1$$
where $V(\bold y)$ is a real-valued potential that is square integrable
on $\Omega,$ and it is assumed that
$V(\bold y)\equiv 0$ outside $\Omega.$
In the spherically-symmetric case in ${\bold R}^3,$ using $V(y)$ instead of
$V(\bold y)$ with $y:=|\bold y|,$ we define the
{\it special transmission eigenvalues}
of (5.1) as
those transmission eigenvalues
for which the
corresponding wavefunctions are spherically symmetric
in addition to $V$ being spherically symmetric.

We remark that the potential $V(y)$ we use in this section does not
necessarily come from any function $\rho$ appearing in (1.1) or (2.1)
via the Liouville transformation in (2.8).
The only assumption we make on $V(y)$ is that it is real valued,
compactly supported within the interval $y\in[0,a],$
and integrable on $(0,a).$ We will say that $V(y)$ belongs to the admissible class $\tilde A$ if $V(y)$ satisfies those conditions.

The results given in the following proposition are either known or can easily be proved
by using the available results [3,6,12,14,18] for the half-line
Schr\"odinger
equation by exploiting the compact-support property of the potential.
We provide a brief proof for the convenience of the reader.

\noindent {\bf Proposition 5.1} {\it
Assume that the potential $V(y)$ belongs to the admissible class
$\tilde A$ with support within the interval $[0,a].$ We then have the following:}

\item{(a)} {\it The corresponding Jost solution
$\tilde f(y;k)$ has an analytic extension from
$k\in{\bold R}$ to the entire complex plane ${\bold C}$ for each fixed $y.$
Similar to (2.11), we have}
$$\tilde f(y;-k)=\tilde f(y;k)^*,\qquad k\in{\bold R}.\tag 5.2$$

\item{(b)} {\it The quantity $\tilde f(0;k)$ is nonzero in
${\bold C^+}$ except perhaps at a finite number of points on the positive
imaginary axis,
say at $k=i\beta_j$ for $j=1,\dots,N$ for some
nonnegative integer $N.$ Such zeros are all simple and
they correspond to the bound states of the half-line Schr\"odinger equation with the Dirichlet
boundary condition at the origin.}

\item{(c)} {\it The large-$k$ asymptotics
of $\tilde f(0;k)$ in ${\bold C}$ are obtained via}
$$\tilde f(0;k)=1+O\left(\displaystyle\frac{1}{k}\right),\qquad k\to \infty\
{\text{in}} \ \overline{\bold C^+}.\tag 5.3$$
$$\tilde f(0;-k)=1+O\left(\displaystyle\frac{1}{k}\right)+e^{-2ika}\,o\left(\displaystyle\frac 1k\right)
,\qquad k\to \infty\
{\text{in}} \ \overline{\bold C^+}.\tag 5.4$$
{\it Hence, $e^{2ikb}\tilde f(0;-k)$ is bounded
in $\overline{\bold C^+}$ for any $b$ satisfying $b\ge a.$}

\item{(d)} {\it The scattering matrix $\tilde S(k)$
for the half-line Schr\"odinger equation with the Dirichlet
boundary condition is defined as [3,6,12,14,18]
$$\tilde S(k):=
\displaystyle\frac{\tilde f(0;-k)}{\tilde f(0;k)},\qquad k\in{\bold R},\tag 5.5$$
and it has a meromorphic extension from
$k\in{\bold R}$ to $k\in{\bold C^+}$ with simple
poles occurring at $k=i\beta_j$ for $j=1,\dots,N.$}

\item{(e)} {\it The scattering matrix satisfies}
$$\tilde S(k)=1+O\left(\displaystyle\frac{1}{k}\right)+e^{-2ika}\,o\left(\displaystyle\frac 1k\right)
,\qquad k\to \infty\
{\text{in}} \ \overline{\bold C^+}.\tag 5.6$$

\item{(f)} {\it
Associated with
each bound state there is a positive number, known as the
corresponding norming constant, defined as}
$$c_j:=\displaystyle\frac{1}{\sqrt{\int_0^\infty dy\,[\tilde f(y;i\beta_j)]^2}},
\qquad j=1,\dots,N.$$
{\it Because $V(y)$ has support confined to the finite interval
$y\in[0,a],$ the norming constants are uniquely determined by $\tilde f(0;k)$ alone,
or equivalently by the scattering matrix $\tilde S(k)$ alone, as
$$c_j=\sqrt{i\,\text{\rm{Res}}\left(\tilde S(k),i\beta_j\right)}
,\qquad j=1,\dots,N,\tag 5.7$$
where $\text{\rm{Res}}\,(\tilde S(k),i\beta_j)$ denotes the residue
of $\tilde S(k)$ at the pole $k=i\beta_j.$}

\item{(g)} {\it The potential $V(y)$ and the Jost solution $\tilde f(y;k)$
are reconstructed from the solution $K(y,\xi)$ to the
Marchenko integral equation}
$$K(y,\xi)+M(y+\xi)+\int_y^\infty ds\,K(y,s)\,M(s+\xi),\qquad
0<y<\xi<+\infty,\tag 5.8$$
{\it where}
$$M(\xi):=\displaystyle\frac{1}{2\pi}\int_{-\infty}^\infty dk\,[1-
\tilde S(k)]\,e^{ik\xi}
+\displaystyle\sum_{j=1}^N c_j^2 e^{-\beta_j\xi}
,\qquad \xi\in{\bold R}.\tag 5.9$$
{\it In fact, we have}
$$V(y)=-2\displaystyle\frac{ dK(y,y)}{dy},\qquad y>0,\tag 5.10$$
$$\tilde f(y;k)=e^{iky}+\int_y^\infty ds\,K(y,s)\,e^{iks}.\tag 5.11$$

\noindent PROOF: The results are mainly known [3,6,12,14,18]; for example,
(5.3) and (5.4) can be
derived by using
the integral representation [3,6,12,14,18] for the Jost solution,
namely by using
$$\tilde f(y;k)=e^{iky}+\displaystyle\frac{1}{k}\displaystyle\int_x^a ds\,\left[\sin k(s-y)\right]\,V(s)\,\tilde f(s;k).$$
We then get (5.6) by using (5.3) and (5.4) in (5.5). The proof of (5.7)
can be outlined as follows. From (2.7) and (5.11), it follows that
$K(y,\xi)=0$ for $a<y<\xi<+\infty,$ and hence (5.8) in turn implies that
$M(y+\xi)=0$ for $a<y<\xi<+\infty.$ This fact, combined with (d) and (e)
allows us to evaluate $M(y+\xi)$ for $a<y<\xi<+\infty$
by using (5.9) as a contour integral along
a semicircle in ${\bold C^+}$ with its
center at the origin and with its
radius becoming infinite in the limit.
Using $M(y+\xi)=0$ for $a<y<\xi<+\infty$ in (5.9), we get (5.7). \qed

Analogously to (1.2), let us define
$$\tilde D(k):=\displaystyle\frac{\sin (ka)}{k}\,\tilde\phi'(a;k)-\cos(ka)\,\tilde\phi(a;k),\tag 5.12$$
where $a$ is the positive constant related to the
support of $V(y)$ and
$\tilde\phi(y;k)$ is the unique solution to (2.16).
The following fundamental result
is the analog of Theorem~2.3, and its proof is omitted because
it is similar to the proof of Theorem~2.3.

\noindent {\bf Theorem 5.2} {\it
Assume that the potential $V(y)$ belongs to the admissible class
$\tilde A.$ Then,
the quantity $\tilde D(k)$ defined in (5.12),
is related to the Jost solution
$\tilde f(y;k)$ appearing in (2.6) and (2.7) as}
$$\tilde D(k)=\displaystyle\frac{\tilde f(0;k)-\tilde f(0;-k)}{2ik},\qquad k\in{\bold C}.\tag 5.13$$
{\it For real $k$-values, we then have}
$$\text{Im}[\tilde f(0;k)]=k\,\tilde D(k),\qquad k\in{\bold R}.\tag 5.14$$

We have [2] the following analog of
Theorem~1.1.

\noindent {\bf Theorem 5.3} {\it Consider the
special case of (5.1) with
$\Omega$ being the three-dimensional ball of
radius $a$ centered at the origin,
where only
spherically-symmetric wavefunctions are allowed
and it is assumed that such wavefunctions
are continuous in the closure of $\Omega.$ Then, the
corresponding special transmission eigenvalues of
(5.1) coincide with the $k^2$-values related to the zeros of the quantity
$\tilde D(k)$ defined in (5.12), where
$\tilde \phi(y;k)$ is the unique solution to
(2.16) with the potential $V(y)$ belonging to the
admissible class $\tilde A.$}

When $V(y)$ belongs to
the admissible class $\tilde A,$ the quantity
$\tilde D(k)$ defined in (5.12) is known [2] to be entire in
$k^2$ and has a
representation analogous to (1.4), namely
$$\tilde D(k)=\tilde\gamma\,
k^{2\tilde d}\prod_{n=1}^{\infty }\left( 1-\frac{
k^2}{\tilde k_n^2}\right),\tag 5.15$$
with $\tilde k_n^2$ for $n\in\bold N$
being the nonzero transmission eigenvalues,
some of which may be repeated, and
$\tilde d$ denoting the multiplicity
of the zero transmission eigenvalue.
As in [2], we refer to the multiplicity
of a nonzero zero $\tilde k_n$ of $\tilde D(k)$ as the
multiplicity of the special transmission eigenvalue
$\tilde k_n^2.$

The following uniqueness result was proved
in [2] and is the analog of Theorem~1.2.

\noindent {\bf Theorem 5.4} {\it Assume that
$V(y)$ belongs to the admissible class $\tilde A.$
Then, $V(y)$ is uniquely determined
by the function $\tilde D(k)$ appearing in (5.12) and (5.15)
if we assume that there exists at least one $V(y)$ in
$\tilde A$ corresponding to
$\tilde D(k).$ Equivalently stated, if the existence
is ensured, $V(y)$ is uniquely determined by the knowledge
of the special transmission eigenvalues of (5.1)
with their multiplicities and the constant $\tilde \gamma$
appearing in (5.15).}

Our goal in this section is to give an independent proof of
Theorem~5.4 and further provide a reconstruction of $V(y)$ from $\tilde D(k).$
The reconstruction consists of the following steps
and the uniqueness follows as a result of the uniqueness in each
reconstruction step.

\item{(a)} First, reconstruct $\tilde f(0;k)$ from $\tilde D(k),$
where $\tilde f(y;k)$ is the Jost solution appearing in (2.6) and (2.7).
This is done by solving the Riemann-Hilbert
problem given by
$$[\tilde f(0;k)-1]-[\tilde f(0;-k)-1]=2ik \tilde D(k),\qquad k\in{\bold R},\tag 5.16$$
which is obtained from (5.13). It follows from (5.3) and (5.16) that $2ik \tilde D(k)$
behaves as $O(1/k)$ as $k\to\pm\infty$ on ${\bold R}.$
Furthermore, because $\tilde D(k)$
has an analytic extension to the entire
complex plane, it follows that $2ik \tilde D(k)$
satisfies the Lipschitz continuity in ${\bold R}.$
Thus, the Riemann-Hilbert problem in (5.16) has
a unique solution that is given by
$$\tilde f(0;k)=1+\displaystyle\frac{1}{\pi}\int_{-\infty}^\infty dt\,
\displaystyle\frac{t \tilde D(t)}{t-k-i0^+},\qquad
k\in \overline{\bold C^+}.\tag 5.17$$
We can write (5.14) as
$$\text{Im}[\tilde f(0;k)-1]=k \tilde D(k),\qquad k\in{\bold R}.\tag 5.18$$
The result given in (5.17) also follows by
using the Schwarz integral formula (4.9), by replacing $Q(k)$ there with
$\tilde f(0;k)-1,$ with the help of (5.18), we obtain (5.17).

\item{(b)} Having obtained $\tilde f(0;k)$ from $\tilde D(k),$ we can
use (5.5) and
(5.7) in (5.9) and obtain the Marchenko kernel
$M(\xi)$ from $\tilde f(0;k).$

\item{(c)} The potential $V(y)$ is then uniquely reconstructed as in (5.10) by using
$M(\xi)$ as input to the Marchenko integral equation (5.8) and by
obtaining $K(y,\xi)$ as the unique solution to the Marchenko equation.

Let us mention that it is an open problem whether the value of
$\tilde\gamma$ appearing in (5.3) can be determined from the
zeros of $\tilde D(k).$ If the answer is yes, then $\tilde\gamma$
is not needed for the unique
determination of $V(y),$ and
the zeros of $\tilde D(k)$ with their multiplicities would be
sufficient for the reconstruction of $V(y).$
In the following example, we show that $\tilde \gamma$
is needed to construct a
potential, which is, however, outside the admissible class $\tilde A.$

\noindent {\bf Example 5.5} Let the potential $V(y)$ be given as
$$V(y)=c\,\delta(y-a),$$
where $c$ is a real nonzero constant, $\delta(y-a)$ denotes the
 Dirac delta function with argument $y-a,$ and $a$ is the positive number related
to the interval $[0,a]$ containing the support of $V(y).$ The corresponding Jost solution is obtained by solving
(2.6) and (2.7), and we get
$$\tilde f(y;k)=\cases \left(1-\displaystyle\frac{c}{2ik}\right)e^{iky}+
\displaystyle\frac{c}{2ik}\,e^{2ika-iky},\qquad y\le a,\\
\noalign{\medskip}
e^{iky},\qquad y\ge a.\endcases\tag 5.19$$
 From (5.19) we evaluate
$\tilde f(0;k)$ and then using (5.13) and (5.15), we obtain
the values of $\tilde \gamma$ and $\tilde D(k),$ yielding
$$\tilde \gamma=ca^2,\quad \tilde
E(k):=\displaystyle\frac{\tilde D(k)}{\tilde \gamma}=\left(\displaystyle\frac{\sin(ka)}{ka}
\right)^2=\displaystyle\prod_{n=1}^\infty \left(1-\displaystyle\frac{a^2k^2}{n^2\pi^2}\right)^2.\tag 5.20$$
Hence, in this example the transmission eigenvalues, i.e. the $k^2$-values
corresponding to the zeros of $\tilde D(k),$
all have double multiplicities and
are given by $k_n^2=n^2\pi^2/a^2$ for $n\in\bold N.$ However, as seen from (5.20)
$\tilde E(k)$ alone does not uniquely determine $c,$ and hence
$c$ or equivalently
$\tilde \gamma$ is also needed for the unique determination of
$V.$

\vskip 8 pt
\noindent {\bf 6. EXAMPLES}
\vskip 3 pt

In this section we illustrate the transmission eigenvalue problem corresponding to $\rho(x)$ appearing in
(1.3) with some explicit examples. In our first example, with the help of (2.2), (2.5), (2.9), and Example~5.5 we present
a concrete $\rho(x)$ for which we can explicitly
evaluate the relevant quantities $D(k)$ and $E(k),$
given in (1.2) and (1.4), respectively.

\noindent {\bf Example 6.1}
Let $\epsilon$ be a positive parameter and let $c$ be a real nonzero
parameter. Assume (2.2) is given by
$$y(x)=\cases \displaystyle\frac{\epsilon^2 x}{\epsilon c x+1},\qquad x\le x_0,\\
\noalign{\medskip}
x-x_0+y_0,\qquad x\ge x_0,\endcases$$
where
$$x_0:=\displaystyle\frac{\epsilon-1}{\epsilon c},\quad
y_0:=y(x_0)=\displaystyle\frac{\epsilon-1}{c}.\tag 6.1$$
In order to have $x_0$ and $y_0$ positive, we must have $c>0$ if
$\epsilon>1$ and we must have $c<0$ if $\epsilon<1.$ Note that (6.1) implies that
$$b-a=x_0-y_0=-\displaystyle\frac{(\epsilon-1)^2}{\epsilon c},\tag 6.2$$
where $a$ and $b$ are the parameters appearing in (1.5). We have $b\ge x_0$ and hence from
(6.2) we see that $a<b$ if $\epsilon>1$ and that $a>b$ if $\epsilon<1.$
Using (2.2) we get
$$\rho(x)=\left(\displaystyle\frac{dy}{dx}\right)^2=
\cases \displaystyle\frac{\epsilon^4}{(\epsilon c x+1)^4},\qquad x\le x_0,\\
\noalign{\medskip}
1,\qquad x\ge x_0,\endcases\tag 6.3$$
Because of (6.1), we see from (6.3) that $\rho(x)$ is continuous
at $x_0,$ whereas $\rho'(x)$ jumps from $\rho'(x_0^-)=-4c$ to
$\rho'(x_0^+)=0.$
One can directly verify that
the Jost solution $f(x;k)$ to (2.1) is given by
$$f(x;k)=\cases
\displaystyle\frac{\epsilon c x+1}{2k\epsilon}\, e^{-ik(\epsilon-1)^2/(\epsilon c)}
Z(x;k,\epsilon,c),\qquad x\le x_0,
\\
\noalign{\medskip}
e^{ikx},\qquad x\ge x_0,\endcases\tag 6.4$$
where we have defined
$$Z(x;k,\epsilon,c):=
(2k+ic)
e^{ik\epsilon^2 x/(\epsilon c x+1)}-ic\,
e^{2ik(\epsilon-1)/c-ik\epsilon^2 x/(\epsilon c x+1)}.$$
One can check that $f(x;k)$ and $f'(x;k)$ are continuous at $x=x_0.$
 From (6.4) we get
$$f(0;k)=\displaystyle\frac{1}{2\epsilon k}\,e^{-ik(\epsilon-1)^2/(\epsilon c)}
\left[(2k+ic)-ic\,e^{2ik(\epsilon-1)/c}
\right].\tag 6.5$$
Using (6.5) in (2.17) we obtain
$$D(k)=\displaystyle\frac{c}{2\epsilon k^2}\left[
\cos\left(\displaystyle\frac{k(\epsilon-1)^2}{\epsilon c} \right)-
\cos\left(\displaystyle\frac{k(\epsilon^2-1)}{\epsilon c} \right)-
\displaystyle\frac{2k}{c}\,\sin\left(\displaystyle\frac{k(\epsilon-1)^2}{\epsilon c} \right)
\right].\tag 6.6$$
By expanding (6.6) in powers of $k^2$ we get
$$D(k)=\gamma k^2+O(k^4),\qquad k\to 0 \ \text {in} \ {\bold C}.$$
where a comparison with (1.4) reveals that
$$d=1,\quad \gamma=-\displaystyle\frac{(\epsilon-1)^4}{3\epsilon^3 c^3}.\tag 6.7$$
Using (6.6) and (6.7) in (1.4) we have
$$E(k)=-\displaystyle\frac{3\epsilon^2 c^4}{2(\epsilon-1)^4 k^2}\left[
\cos\left(\displaystyle\frac{k(\epsilon-1)^2}{\epsilon c} \right)-
\cos\left(\displaystyle\frac{k(\epsilon^2-1)}{\epsilon c} \right)-
\displaystyle\frac{2k}{c}\,\sin\left(\displaystyle\frac{k(\epsilon-1)^2}{\epsilon c} \right)
\right].\tag 6.8$$
 From (6.8), we get
$$E(k)=\displaystyle\frac{3\epsilon^2 c^3}{(\epsilon-1)^4k}\sin \left(\displaystyle\frac{k(\epsilon-1)^2}{\epsilon c} \right)+O\left(\displaystyle\frac{1}{k^2}\right),\qquad k\to\pm\infty,\tag 6.9$$
and hence a comparison of (6.9) with (3.2) reveals that
$$b-a=-\displaystyle\frac{(\epsilon-1)^2}{\epsilon c},\quad
\gamma\,[\rho(0)]^{1/4}=-\displaystyle\frac{(\epsilon-1)^4}{3\epsilon^2 c^3},$$
which is compatible with
the value of $(b-a)$ given in (6.2),
$\rho(0)$ from (6.3), and $\gamma$ in (6.7).
Note that $\rho(x)$ in this example is outside the admissible class $\Cal A$ because
of the jump discontinuity of $\rho'(x)$ at $x_0.$
Let us remark that $D(k)$ given in (6.6) can also be obtained by using (1.2),
where
the unique solution
$\phi(x;k)$ to (2.15) in this case is given by
$$\phi(x;k)=\cases
\displaystyle\frac{ (\epsilon c x+1)}{\epsilon^2 k}\,\sin\left(\displaystyle\frac{\epsilon^2 kx}
{\epsilon c x+1}\right),\qquad x\le x_0,
\\
\noalign{\medskip}
c_3(k,\epsilon,c)\,\sin(kx)+c_4(k,\epsilon,c)\,\cos(kx),
\qquad x\ge x_0
,\endcases$$
with the constants $c_3(k,\epsilon,c)$ and $c_4(k,\epsilon,c)$ specified as
$$c_3(k,\epsilon,c):=\displaystyle\frac{1}{2\epsilon k}\,\cos\left(
\displaystyle\frac{k(1-\epsilon)^2}{\epsilon c}\right)
+\displaystyle\frac{c}{4\epsilon k^2}
\left[\sin\left(
\displaystyle\frac{k(1-\epsilon)^2}{\epsilon c}\right)
-\sin\left(
\displaystyle\frac{k(1-\epsilon^2)}{\epsilon c}\right)
\right],$$
$$c_4(k,\epsilon,c):=\displaystyle\frac{1}{2\epsilon k}\,\sin\left(
\displaystyle\frac{k(1-\epsilon)^2}{\epsilon c}\right)
-\displaystyle\frac{c}{4\epsilon k^2}
\left[\cos\left(
\displaystyle\frac{k(1-\epsilon)^2}{\epsilon c}\right)
-\cos\left(
\displaystyle\frac{k(1-\epsilon^2)}{\epsilon c}\right)
\right].$$

Using the result of Example~6.1, in the next example we will produce two distinct
profiles $\rho(x)$ corresponding to the same $E(k)$ but to different $\gamma$ values;
in fact, in one case we will have $a>b$ and in the other case we will have $a<b.$

\noindent {\bf Example 6.2}
In Example~6.1 above, let us use the following values for the parameters
$$\epsilon=2, \quad c=\displaystyle\frac{1}{b},$$
where $b$ is the constant that appears in (1.2)
and is related to the known support of $\rho(x)-1.$
Using (6.2), (6.3), (6.7), and (6.8), we obtain
$$\rho(x)=
\cases \displaystyle\frac{16b^4}{(2x+b)^4},\qquad x\le b/2,\\
\noalign{\medskip}
1,\qquad x\ge b/2,\endcases\tag 6.10$$
$$\gamma=\displaystyle\frac{-b^3}{24},\quad a=\displaystyle\frac{3b}{2},\tag 6.11$$
$$E(k)=\displaystyle\frac{-6}{b^4 k^2}\left[
\cos\left(\displaystyle\frac{bk}{2}\right)-\cos\left(\displaystyle\frac{3bk}{2}\right)-2bk\,\sin\left(\displaystyle\frac{bk}{2}\right)
\right].\tag 6.12$$
On the other hand, in Example~6.1 if we use the parameters
$$\epsilon=\displaystyle\frac{1}{2}, \quad c=-\displaystyle\frac{1}{b},$$
 from (6.2), (6.3), (6.7), and (6.8), then we obtain
$$\rho(x)=
\cases \displaystyle\frac{b^4}{(2b-x)^4},\qquad x\le b,\\
\noalign{\medskip}
1,\qquad x\ge b,\endcases\tag 6.13$$
$$\gamma=\displaystyle\frac{b^3}{6},\quad a=\displaystyle\frac{b}{2},\tag 6.14$$
$$E(k)=\displaystyle\frac{-6}{b^4 k^2}\left[
\cos\left(\displaystyle\frac{bk}{2}\right)-\cos\left(\displaystyle\frac{3bk}{2}\right)-2bk\,\sin\left(\displaystyle\frac{bk}{2}\right)
\right].\tag 6.15$$
Thus, as seen from (6.12) and (6.15), we have produced two distinct profiles for $\rho(x)$ given in (6.10) and (6.13), respectively, corresponding to the same $E(k),$ but two different $\gamma$ values.
In fact, as seen from (6.11) and (6.14), the former corresponds to the case
$a>b$ and the latter to $a<b.$ We can simplify and rewrite (6.12) as
$$E(k)=\displaystyle\frac{12}{b^3k}\,\sin\left(\displaystyle\frac{bk}{2}\right)
\left(1-\displaystyle\frac{\sin(bk)}{bk}\right).\tag 6.16$$
As seen from (6.16), corresponding to the two
distinct profiles given in (6.10) and (6.13), we have
a simple zero transmission eigenvalue, infinitely many simple nonzero real
transmission eigenvalues that are given by
$k_n^2=4n^2\pi^2/b^2$ for $n\in\bold N,$ and infinitely many simple complex
transmission eigenvalues that are related to nonzero zeros of $kb-\sin(kb).$
Note that for each complex transmission eigenvalue, its complex conjugate is
also a transmission eigenvalue.

We conclude with another explicit example.

\noindent {\bf Example 6.3} For a positive parameter $c,$ let
$$\rho(x)=
\cases \displaystyle\frac{(b+c)^2}{(x+c)^2},\qquad x\le b,\\
\noalign{\medskip}
1,\qquad x\ge b,\endcases$$
where $b$ is the positive parameter appearing in (1.5).
Using (2.2) we obtain
$$y(x)=\cases (b+c)\,\log\left(1+\displaystyle\frac{x}{c}\right),\qquad x\le b,\\
\noalign{\medskip}
x-b+a,\qquad x\ge b,\endcases\tag 6.17$$
where $a$ is the parameter appearing in (1.5) and its value
is obtained from (6.17) as
$$a:=y(b)=(b+c)\,\log\left(1+\displaystyle\frac{b}{c}\right).\tag 6.18$$
In this case $a>b$ because
 from (6.18) it follows that
$$\displaystyle \frac{a}{b}=\left(1+\displaystyle\frac{c}{b}\right)\log\left(1+\displaystyle\frac{b}{c}\right)>1,$$
based on the positivity assumption on $c.$
We can solve (1.3) explicitly and get
$$\phi(x;k)=\displaystyle\frac{1}{r_+-r_-}\left[-c^{r_+}(x+c)^{r_-}+c^{r_-}(x+c)^{r_+}\right],\qquad
x\in[0,b],\tag 6.19$$
where we have defined
$$r_{\pm}:=\displaystyle\frac{1}{2}\left[1\pm\displaystyle\sqrt{1-4(b+c)^2 k^2}\right].\tag 6.20$$
Using (6.19) in (1.2) we obtain
$$D(k)=\displaystyle\frac{\sin(bk)}{k}\,\phi'(b;k)-\cos(kb)\,\phi(b;k),\tag 6.21$$
where as we see from (6.19) and (6.20), with the help
of $r_++r_-=1,$ we have
$$\phi(b;k)=\displaystyle\frac{c}{r_+-r_-}\left[-\left(1+\displaystyle\frac{b}{c}\right)^{r_-}
+\left(1+\displaystyle\frac{b}{c}\right)^{r_+}\right],$$
$$\phi'(b;k)=\displaystyle\frac{c}{r_+-r_-}\left[-\displaystyle\frac{r_-}{b+c}\left(1+\displaystyle\frac{b}{c}\right)^{r_-}
+\displaystyle\frac{r_+}{b+c}\left(1+\displaystyle\frac{b}{c}\right)^{r_+}\right].$$
Letting $k\to 0$ in (6.21), with the help of (1.4), we get
$$d=1,\quad \gamma=c^3\left[
-\displaystyle\frac23\,\left(\displaystyle\frac{b}{c}\right)^3-3\left(\displaystyle\frac{b}{c}\right)^2-2\left(\displaystyle\frac{b}{c}\right)
+2\left(1+\displaystyle\frac{b}{c}\right)^2\log\left(1+\displaystyle\frac{b}{c}\right)\right],\tag 6.22$$
where, by using a graphical argument, it can be shown that $\gamma<0.$
We remark that the results in this example are also valid if $c<0$ but $b>-c.$ In that case from (6.18) we get $a<b$ and from (6.22) we get $\gamma>0.$

\vskip 8 pt

\noindent {\bf Acknowledgments.} The first author has been partially supported by
DOD-BC063989 and he is grateful for the hospitality he received
during a recent visit to the National Technical University of Athens.
The second author has been partially supported by a $\Pi$.E.B.E. grant from the National Technical University of Athens. The second author is grateful to the
Washington University in St. Louis, where he is currently a Visiting Professor/Research Fellow
at the Boeing Center for Technology, Information and Manufacturing in the
Olin School of Business.

\vskip 5 pt

\noindent {\bf{References}}


\item{[1]}
L. V. Ahlfors,
{\it Complex analysis,} 3rd ed., McGraw-Hill, New York, 1979.

\item{[2]}
T. Aktosun, D. Gintides, and V. G. Papanicolaou, {\it  The uniqueness in the inverse problem for transmission eigenvalues for the spherically symmetric variable-speed wave equation,}
Inverse Problems {\bf 27}, 115004 (2011).

\item{[3]} T. Aktosun and R. Weder, {\it
Inverse spectral-scattering problem with two sets of discrete spectra for the radial Schr\"odinger equation,} Inverse Problems {\bf 22}, 89--114
(2006).

\item{[4]}
F. Cakoni, D. Colton, and H. Haddar,
{\it On the determination of Dirichlet or transmission eigenvalues from far field data,}
C. R. Math. Acad. Sci. Paris {\bf 348}, 379--383 (2010).

\item{[5]}
F. Cakoni, D. Colton, and P. Monk,
{\it
On the use of transmission eigenvalues to estimate the index of refraction from far field data,} Inverse Problems {\bf 23}, 507--522 (2007).

\item{[6]} K. Chadan and P. C. Sabatier, {\it Inverse problems in
quantum scattering theory,} 2nd ed., Springer, New York, 1989.

\item{[7]} E. A. Coddington and N. Levinson, {\it Theory of ordinary
differential equations,} McGraw-Hill, New York, 1955.

\item{[8]}
D. Colton and R. Kress,
{\it Inverse acoustic and electromagnetic scattering theory,}
2nd ed.,  Springer,
New York, 1998.

\item{[9]}
D. Colton and P. Monk,
{\it The inverse scattering problem for time-harmonic acoustic waves in an inhomogeneous medium,}
Quart. J. Mech. Appl. Math. {\bf 41}, 97--125 (1988).

\item{[10]}
D. Colton, L. P\"aiv\"arinta, and J. Sylvester,
{\it The interior transmission problem,}
Inverse Probl. Imaging {\bf 1}, 13--28 (2007).

\item{[11]} F. D. Gakhov, {\it Boundary value problems,}
Pergamon Press, Oxford, 1966.

\item{[12]} V. A. Marchenko,
{\it Sturm-Liouville operators and applications,}
Birkh\"auser, Basel, 1986.

\item{[13]} N. I. Muskhelishvili, {\it Singular integral equations,}
Wolters-Noordhoff Publishing, Groningen, the Netherlands, 1958.

\item{[14]} B. M. Levitan, {\it
Inverse Sturm Liouville Problems,} VNU Science Press, Utrecht, 1987.

\item{[15]}
J. R. McLaughlin and P. L. Polyakov,
{\it On the
uniqueness of a spherically symmetric speed of sound from transmission
eigenvalues,} J. Differential Equations {\bf 107}, 351--382 (1994).

\item{[16]}
J. R. McLaughlin, P. L. Polyakov, and P. E. Sacks,
{\it Reconstruction of a spherically symmetric speed of sound,}
SIAM J. Appl. Math. {\bf 54}, 1203--1223 (1994).

\item{[17]}
J. R. McLaughlin, P. E. Sacks, and M. Somasundaram,
{\it Inverse scattering in acoustic media using interior transmission eigenvalues,}
in: G. Chavent, G. Papanicolaou, P. Sacks, and W. Symes (eds.),
{\it Inverse problems in wave propagation,} Springer,
New York, 1997, pp. 357--374.

\item{[18]} R. G. Newton, {\it
Scattering theory of waves and particles,} 2nd ed., Springer, New York, 1982.

\end